\documentclass[12pt,reqno,fleqn]{amsart}
\usepackage{amssymb}
\usepackage{amsxtra}

\theoremstyle{plain}
\newtheorem{Thm}{Theorem}

\newtheorem{Prop}{Proposition}
\theoremstyle{definition}
\newtheorem{Def}{Definition}
\theoremstyle{remark}
\newtheorem{Rem}{Remark}

\def\eop {\hfill $\square$}

\numberwithin{equation}{section}

\newcommand{\R}{\mathbf{R}}

\newcommand{\sinc}{\mathrm{sinc}}
\newcommand{\sign}{\mathrm{sign}}
\newcommand{\Vol}{\mathrm{Vol}}

\begin{document}

\title[Sinc Integrals]{Some Remarks on Sinc Integrals
and their Connection with Combinatorics, Geometry and Probability}

\date{January 7, 2002}

\author{David~M. Bradley}
\address{Department of Mathematics \& Statistics\\
         University of Maine\\
         5752 Neville Hall
         Orono, Maine 04469-5752\\
         U.S.A.}
\email{dbradley@member.ams.org, bradley@math.umaine.edu}
\urladdr{http://www.math.umaine.edu/faculty/bradley/}
\thanks{Published in \textit{Analysis (Munich)} \textbf{22} (2002), no.~2, 219--224. [MR
1916426] (2003i:42012)}

\subjclass{Primary: 42A38; Secondary: 52B11, 42A61, 42A85, 26B15,
51M25, 60E10}

\keywords{Sinc integrals, sine integrals, Fourier transform,
characteristic function, Fourier inversion, convolution,
probability density, distribution of a sum, cumulative
distribution, inclusion-exclusion principle.}

\begin{abstract}
   We give an alternative, combinatorial/geometrical evaluation of
   a class of improper sinc integrals studied by the Borweins.  A
   probabilistic interpretation is also noted and used to shed light
   on a related combinatorial identity.
\end{abstract}

\maketitle

\section{Introduction}\label{sect:Intro}
In~\cite{BB}, among other things the Borweins study integrals of
the form
\begin{equation}
   \qquad I_n:=\int_0^\infty \prod_{j=0}^n\sinc(a_jx)\,dx,
\label{BorSinc}
\end{equation}
where $\sinc(x):=(\sin x)/x$ for $x\ne 0$ and $\sinc(0):=1$. They
use a version of the Parseval/Plancherel formula to prove that,
subject to certain conditions on the parameters $a_j$, the
sequence $I_1,I_2,\dots$ is decreasing;  they also give an
explicit evaluation of $I_n$ by expanding the product of sines
into a sum of cosines and then integrating by parts.  As the
Borweins observed, the integral~(\ref{BorSinc}) can be interpreted
geometrically as the volume of an $n$-dimensional polyhedron
obtained by cutting an $n$-dimensional hypercube by two
$(n-1)$-dimensional hyperplanes, and it is interesting to note
that their formula is essentially a signed sum over the vertices
of the hypercube.  It may therefore be of interest to give an
independent evaluation of~(\ref{BorSinc}) using methods from
combinatorial geometry.  The derivation is provided in
\S\ref{sect:Volume}.

In~\cite{BradGupta}, related integrals arise in a probabilistic
setting---namely the problem of determining the probability
density of a sum of independent random variables uniformly
distributed in different ranges.  We employ this interpretation
here to give an alternative proof of an elegant combinatorial
identity related to the evaluation of $I_n$ using methods from
probability theory.  See Proposition~\ref{prop:CoolIdentity}
below.

\section{Fourier Background and Connection with Probability Theory}
\label{sect:Background}
For our purposes, it is more convenient to take products from $1$
to $n$, as opposed to $0$ to $n$ in~(\ref{BorSinc}).  Thus, we fix
a positive integer $n$, and let $\vec a=(a_1,a_2,\dots,a_n)$ be a
vector of positive real numbers. For each $j=1,2,\dots,n$, define
a step function $\chi_j:\R\to\R$ by
\begin{equation}
   \qquad 2a_j\,\chi_j(x) = \left\{\begin{array}{lll}1 &\mbox{if
   $|x|<a_j$,}\\ \tfrac12 &\mbox{if $|x|=a_j$,}\\ 0 &\mbox{if
   $|x|>a_j$.}\end{array}\right.
\label{ChiDef}
\end{equation}
Let $f_n:\R\to\R$ denote the $n$-fold convolution
$\chi_1*\chi_2*\dots*\chi_n$, so that for all real $x$,
\[
  \qquad f_n(x) =
  \int_{-\infty}^\infty\chi_1(x-y_2)\int_{-\infty}^\infty
   \chi_2(y_2-y_3)
   \cdots\int_{-\infty}^\infty\chi_{n-1}(y_{n-1}-y_n)
   \chi_n(y_n)\,dy_2\dots dy_n.
\]
It is now easy to see that $f_n(x)$ is simply the volume (in the
sense of Lebesgue measure on $\R^{n-1})$ of the
$(n-1)$-dimensional polyhedron
\begin{equation}
   \qquad\{(x_1,x_2,\dots,x_{n})\in\R^{n} :
   \textstyle\sum_{j=1}^{n} x_j=x
   \;{\mathrm{and}}\; |x_j|<a_j\;
   {\mathrm{for}}\; j=1,2,\dots,n\}
\label{polyhedron}
\end{equation}
divided by the volume $\prod_{j=1}^n (2a_j)$ of the
$n$-dimensional hypercube
\begin{equation}
   \qquad
   \{(x_1,x_2,\dots,x_n)\in\R^n : |x_j|<a_j \quad {\mathrm{for}}
   \quad j=1,2,\dots,n\}.
\label{hypercube}
\end{equation}
The volume of the region~(\ref{polyhedron}) will be evaluated
directly in \S\ref{sect:Volume}; but first we connect this problem
with that of evaluating integrals such as~(\ref{BorSinc}). The
Fourier transform of $f_n$ is given by
\begin{equation}
   \qquad \widehat{f}_n(t) = \prod_{j=1}^n \widehat{\chi}_j(t)
   =\prod_{j=1}^n\int_{-\infty}^\infty e^{itx}\chi_j(x)\,dx
   =\prod_{j=1}^n\sinc(a_j t),\qquad t\in\R.
\label{sincprod}
\end{equation}
Since for all real $x$, $f_n(x) = \tfrac12 f_n(x+) + \tfrac12
f_n(x-)$ by continuity of $f_n$ for $n>1$, and by definition of
$\chi_1$ when $n=1$, Fourier inversion gives
\begin{equation}
   \qquad f_n(x) = \frac1{2\pi}\int_{-\infty}^\infty e^{-itx}
   \widehat{f}_n(t)\,dt =
   \frac1{2\pi}\int_{-\infty}^\infty e^{itx}\prod_{j=1}^n
   \sinc(a_j t)\,dt,\qquad x\in\R,
\label{SincIntegral}
\end{equation}
where we have made the change of variable $t\mapsto-t$ and used
the fact that $\sinc$ is even.  We remark that apart from our
convention of starting products at $j=1$, $f_n(0)$ coincides with
the Borweins' integral~(\ref{BorSinc}) up to a numerical factor.
More generally, $f_n(x)$ falls under the scope of their Theorem
2(ii). However, to facilitate the presentation of our
combinatorial evaluation, it is convenient to reformulate their
result as follows.  We first make the following
\begin{Def}\label{TauDef} Let $\tau:\R\to\R$ be given by
\begin{equation}
   \qquad\tau(x) = \left\{\begin{array}{lll} 1,& \mbox{if $x>0$,}\\
   \tfrac12, & \mbox{if $x=0$,}\\ 0, &\mbox{if
   $x<0$,}\end{array}\right.
\label{taudef}
\end{equation}
and for $x$ real and $n$ a positive integer, let $x^{n-1}_+:=
x^{n-1}\tau(x)$.  Note that $x^0_+ = \tau(x)$ and
$x^n_+=(\max(x,0))^n$ for $n>0$.
\end{Def}
We are now ready to state our reformulation.
\begin{Thm}\label{Thm:tau}
Let $n$ be a positive integer, let $\vec a=(a_1,a_2,\dots,a_n)$ be
a vector of positive real numbers, and let $f_n(x)$ be as
in~(\ref{SincIntegral}).  Then for all real $x$,
\begin{equation}
\label{CtsFormula2}
\begin{split}
  \qquad f_n(x) &=
  \bigg[\sum_{\vec\varepsilon\in\{-1,1\}^n}
    (x+\vec\varepsilon\cdot\vec a)^{n-1}_+
  \;\prod_{j=1}^n\varepsilon_j\bigg]\bigg/
  \bigg[(n-1)!\,\prod_{j=1}^n (2a_j)\bigg],
\end{split}
\end{equation}
in which the sum is over all $2^n$ vectors of signs
$\vec\varepsilon=(\varepsilon_1,\varepsilon_2,\dots,\varepsilon_n)$
with $\varepsilon_j=\pm 1$ for each $j=1,2,\dots,n$, and of course
$\vec\varepsilon\cdot\vec a$ denotes the inner product
$\sum_{j=1}^n \varepsilon_j a_j$.
\end{Thm}
To see more directly the connection with Theorem 2(ii)
of~\cite{BB}, simply replace the $\tau$ function in
Theorem~\ref{Thm:tau} by the sign function using the relationship
$\sign(x)=2\tau(x)-1$ for $x$ real.  It then suffices to prove the
identity
\begin{equation}
\label{CoolIdentity}
   \qquad\sum_{\vec\varepsilon\in\{-1,1\}^n}
   (x+\vec\varepsilon\cdot\vec a)^{n-1}\;\prod_{j=1}^n
   \varepsilon_j = 0,\qquad x\in\R,
\end{equation}
which is a special case of~\cite[Theorem 2(i)]{BB}, there proved
most elegantly using the method of generating functions.  We
conclude this section with an alternative proof
of~(\ref{CoolIdentity}) based on probabilistic considerations.
For the reader's convenience, the result is restated as
\begin{Prop}[Theorem 2(i) of~\cite{BB}]\label{prop:CoolIdentity}
Let $n$ be a positive integer and $\vec a=(a_1,a_2,\dots,a_n)$ a
vector of positive real numbers.  Then for all real $x$, we have
\[
   \qquad\sum_{\vec\varepsilon\in\{-1,1\}^n}
   (x+\vec\varepsilon\cdot\vec a)^{r}\;\prod_{j=1}^n
   \varepsilon_j =
   \left\{\begin{array}{ll} 0,&\quad\mbox{if\; $r=0,1,2,\dots,n-1$,}\\
   n!\,2^n\displaystyle\prod_{j=1}^n a_j,&\quad
   \mbox{if\; $r=n$.}\end{array}\right.
\]
\end{Prop}

\begin{proof}
As in~\cite{BradGupta}, we note that $f_n(x)$ represents the
probability density of the sum of $n$ independent random variables
$X_1,X_2,\dots,X_n$ with $X_j$ uniformly distributed in the
interval $[-a_j,a_j]$ for each $j=1,2,\dots,n$.  In light of
Theorem~\ref{Thm:tau}, the cumulative distribution
$F_n(x):=\Pr(\sum_{j=1}^n X_j \le x)$ is given by
\begin{equation}
   \qquad F_n(x) = \int_{-\infty}^x f_n(y)\,dy
   =  \bigg[\sum_{\vec\varepsilon\in\{-1,1\}^n}
  (x+\vec\varepsilon\cdot\vec a)^{n}_+
  \;\prod_{j=1}^n\varepsilon_j\bigg]\bigg/
  \bigg[n!\,\prod_{j=1}^n (2a_j)\bigg].
\label{Cumulative}
\end{equation}
Let $A_n:=\sum_{j=1}^n a_j$. Since each $X_j$ is uniformly
distributed in $[-a_j,a_j]$, the sum $\sum_{j=1}^n X_j$ must fall
within the interval $[-A_n,A_n]$, whence
$F_n(x)=F_n(\infty)=\Pr(\sum_{j=1}^n X_j\le A_n)=1$ for all $x\ge
A_n$.  But if $x\ge A_n$, then $x+\vec\varepsilon\cdot\vec a\ge 0$
for each $\vec\varepsilon\in\{-1,1\}^n$, and so for such $x$ we
can drop the subscripted ``+'' from~(\ref{Cumulative}). It follows
that
\begin{equation}
   \qquad\sum_{\vec\varepsilon\in\{-1,1\}^n}
  (x+\vec\varepsilon\cdot\vec a)^{n}
  \;\prod_{j=1}^n\varepsilon_j
  =  n!\,\prod_{j=1}^n (2a_j)
  =  n!\,2^n\prod_{j=1}^n a_j
\label{ConstantPoly}
\end{equation}
holds for all $x\ge A_n$. Since the left hand side
of~(\ref{ConstantPoly}) is a polynomial of degree at most $n$ in
$x$ which is constant for all sufficiently large values of $x$, it
must in fact be constant for all real $x$ by the identity theorem.
This proves the case $r=n$ of the proposition.  The other cases
follow by repeatedly differentiating with respect to $x$.
\end{proof}

\begin{Rem} A related probabilistic interpretation of the
integral~(\ref{BorSinc}), which is essentially our $f_n(0)$, can be
found in~\cite[p.\ 81, Remarks 1(b)(ii)]{BB}.
\end{Rem}

\section{Combinatorial Proof of Theorem~\ref{Thm:tau}}
\label{sect:Volume}

We will prove the formula~(\ref{CtsFormula2}) of
Theorem~\ref{Thm:tau} by evaluating the volume of the
polyhedron~(\ref{polyhedron}) directly.  Let us denote this volume
by $V_n(x)$.  Since the case $n=1$ is trivial, we will assume
$n>1$. Let the vector $\vec a$ be as in the previous section, and
let $N$ be the set of vectors in $\R^n$ whose respective
components exceed the components of $-\vec a$. Thus,
$N:=\{(x_1,x_2,\dots,x_n)\in\R^n : x_j> -a_j \quad {\mathrm{for}}
\quad j=1,2,\dots,n\}$.  To each subset $A$ of $N$, we associate a
non-negative function $w(A):\R\to\R$ defined by
\[
   \qquad w(A)(x) := \Vol\{(x_1,x_2,\dots,x_n)\in A: \textstyle\sum_{j=1}^n
   x_j=x\},
\]
where, as in~(\ref{polyhedron}), the volume is taken in the sense
of Lebesgue measure on $\R^{n-1}$.  Clearly $w$ is a non-negative
additive weight-function on the subsets of $N$; for if $A$ and $B$
are disjoint subsets of $N$, then $w(A\cup B)=w(A)+w(B)$, and
$w(A)\ge 0$ for all $A\subseteq N$.  Therefore, if we let
$S:=\{1,2,\dots,n\}$ and $R_k:=\{(x_1,x_2,\dots,x_n)\in N: x_k>
a_k\}$ for $k\in S$, then
\[
   \qquad\bigcap_{j\in S}(N\setminus R_j) = \{(x_1,x_2,\dots,x_n)\in \R^n:
   - a_j<x_j\le a_j \quad {\mathrm{for}}\quad j=1,2,\dots,n\},
\]
and hence in view of~(\ref{polyhedron}) and the definition of $w$,
we have for every $x\in\R$,
\begin{equation}
   \qquad V_n(x) = w\bigg[\bigcap_{j\in S}(N\setminus R_j)\bigg](x),
\label{WeightRep}
\end{equation}
since the underlying sets differ by a set of measure zero in
$\R^{n-1}$.  The Inclusion-Exclusion Principle states that
\begin{equation}
   \qquad w\bigg[\bigcap_{j\in S}(N\setminus R_j)\bigg]
   =\sum_{T\subseteq S} (-1)^{|T|} \,w\big(\bigcap_{j\in T}
   R_j\big),
\label{IEP}
\end{equation}
where the sum is over all subsets $T$ of $S$, $|T|$ denotes the
cardinality of $T$, and under the usual convention for
intersections over empty sets, the term corresponding to $T=\{\}$
is $w(N)$.  Fortunately, the weights of the sets $R_k$ and their
intersections are relatively easy to compute. We make a linear
change of variable
\[
   \qquad y_j := \sum_{k=j}^n x_k, \qquad j=1,2,\dots,n.
\]
Then $y_j=y_{j+1}+x_j$ for $j=1,2,\dots,n-1$, $y_n=x_n$, and the
Jacobian of the transformation is 1.  Given $T\subseteq S$, let
$\vec b=\vec b(T) = (b_1,b_2,\dots,b_n)$ be the vector whose $j$th
component $b_j$ is equal to $a_j$ if $j\in T$, and $b_j=-a_j$
otherwise.  Then
\[
   \qquad R(T) := \bigcap_{j\in T}R_j = \{(x_1,x_2,\dots,x_n)\in \R^n :
   x_j > b_j \quad{\mathrm{for}}\quad j=1,2,\dots,n\}.
\]
For each $j=1,2,\dots,n$, let $c_j := \sum_{k=j}^n b_j$.  If $\vec
x\in R(T)$, then $y_j>c_j$, and $y_j-b_j> y_{j+1}>c_{j+1}$ for
$j=1,2,\dots,n-1$.  Thus,
\[
   \qquad w(R(T))(y_1) = \int_{c_2}^{y_1-b_1} \int_{c_3}^{y_2-b_2}
  \cdots\int_{c_n}^{y_{n-1}-b_{n-1}} dy_n \cdots dy_3\, dy_2.
\]
Now make the change of variable $z_j=y_j-c_j$ for $j=1,2,\dots,n$.
Then each $z_j>0$, and we obtain
\[
   \qquad w(R(T))(y_1) =
   \int_0^{z_1}\int_0^{z_2}\cdots\int_0^{z_{n-1}}
   dz_n\cdots dz_3\,dz_2
   = \frac{z_1^{n-1}}{(n-1)!}.
\]
In light of the constraint $y_1-c_1=z_1>0$, we see that
\[
   \qquad w\big(\bigcap_{j\in T}R_j\big)(y_1)
   =w(R(T))(y_1) = \frac{(y_1-c_1)^{n-1}_+}{(n-1)!}
   =\frac{(y_1-\sum_{j=1}^n b_j)^{n-1}_+}{(n-1)!}.
\]
Thus, from~(\ref{WeightRep}) and~(\ref{IEP}) it follows that
\[
   \qquad V_n(x)=\sum_{T\subseteq S} (-1)^{|T|} \,w\big(\bigcap_{j\in T}
   R_j\big)(x)=\frac{1}{(n-1)!}\sum_{\vec\varepsilon\in\{-1,1\}^n}
   (x+\vec\varepsilon\cdot\vec a)^{n-1}_+\;\prod_{j=1}^n\varepsilon_j,
\]
where we have made the correspondence $\varepsilon_j=-1$ if $j\in
T$ and $\varepsilon_j=1$ if $j\notin T$ in each term of the sums.
Since $f_n(x)=V_n(x)/\prod_{j=1}^n (2a_j)$, the proof of
Theorem~\ref{Thm:tau} is complete. \eop

\end{document}